\documentclass[review,onefignum,onetabnum]{article}

\usepackage{graphicx, fancyhdr, color} 
\usepackage{epstopdf}
\usepackage{amssymb}
\usepackage{amsmath}
\usepackage{amsfonts}
\usepackage{subfigure}
\usepackage{enumerate}
\usepackage{algorithmic}
\usepackage[ruled,vlined,algosection]{algorithm2e}
\usepackage{comment}
\usepackage{hyperref}

\newcommand{\norm}[1]{\left\Vert#1\right\Vert}

\usepackage[a4paper, total={6in, 8in}]{geometry}

\newtheorem{theorem}{Theorem}[section]

\newtheorem{lemma}[theorem]{Lemma}

\title{Inexact Newton Method for General Purpose Reservoir Simulation}

\author{Soham Sheth\thanks{Geoscience Research Centre, TOTAL E\&P UK, Westhill, UK. \href{mailto:soham-sheth@utulsa.edu}{soham-sheth@utulsa.edu}, \href{mailto:arthur.moncorge@total.com}{arthur.moncorge@total.com}} \qquad Arthur Moncorg\'e\footnotemark[1]}

\begin{document}
\maketitle

\begin{abstract}
Inexact Newton Methods are widely used to solve systems of nonlinear equations. The convergence of these methods is controlled by the relative linear tolerance, $\eta_\nu$, that is also called the forcing term. A very small $\eta_\nu$ may lead to \textit{oversolving} the Newton equation. Practical reservoir simulation uses inexact Newton methods with fixed forcing term, usually in the order of $10^{-3}$ or $10^{-4}$. Alternatively, variable forcing terms for a given inexact Newton step have proved to be quite successful in reducing the degree of \textit{oversolving} in various practical applications. The cumulative number of linear iterations is usually reduced, but the number of nonlinear iterations is usually increased. 
We first present a review of existing inexact Newton methods with various forcing term estimates and then we propose improved estimates for $\eta_\nu$. 
These improved estimates try to avoid as much as possible \textit{oversolving} when the iterate is far from the solution and try to enforce quadratic convergence in the neighborhood of the solution.
Our estimates reduce the total linear iterations while only resulting in few extra Newton iterations. 
We show successful applications to fully-coupled three-phase and multi-component multiphase models in isothermal and thermal steam reservoir simulation
as well as a real deep offshore west-African field with gas re-injection 
using the reference CPR-AMG iterative linear solver. 
\end{abstract}


\section{Introduction}
\label{Sec:Intro}
Various computational engineering problems require the solution of a system of highly nonlinear functions given by 
\begin{equation}\label{Eq:nonlinear}
    \mathcal{R}(u) = 0
\end{equation}
with the nonlinear operator $\mathcal{R}: D \subset \mathbb{R}^n \rightarrow \mathbb{R}^n$, the domain $D$ in $\mathbb{R}^n$ and the state variables solution $u \in D$.
Given an initial guess $u^0$, Eq.~\eqref{Eq:nonlinear}
can be solved with the Newton method \cite{deuflhard2004,Ortega2000}. 
The Newton method can be written as a sequence of successive linearizations of the form 
\begin{eqnarray} 
    \mathcal R'(u^\nu) \, \delta^\nu &=& - \mathcal R(u^\nu),\label{Eq:General_Newton} \\
    u^{\nu+1} &=& u^\nu + \delta^\nu 
\end{eqnarray}
with $u^{\nu}$, $u^{\nu+1}$ the old and the new iterates of the solution,
$\delta^{\nu}$ the update of the solution and
$\mathcal R'(u^\nu)$ the Jacobian matrix that can be very large and sparse. 
The reader is referred to \cite{PeacemanBook,aziz79} for an introduction on the equations and the models used in reservoir simulation.
The current standard in reservoir simulation is the Generalized Minimal Residual method (GMRES) \cite{saad1986gmres}
preconditioned by CPR\cite{wallis85,Lacroix2003Decoupl,Cao2005CPR}-AMG \cite{stuben2001a, henson2002}. 
We refer a Newton iteration as an outer iteration and we refer a GMRES iteration as an inner iteration. A very tight tolerance for the inner iteration results in the exact solution of the Newton equation given by Eq.~\eqref{Eq:General_Newton}. This may lead to \textit{oversolving} the Newton equation, which means imposing accuracy that leads to significant disagreement between $\mathcal{R}(u^{\nu+1})$ and it's local linear model, $\mathcal R'(u^\nu) \delta^\nu + \mathcal R(u^\nu)$. This may entail little to no decrease in the nonlinear residual norm over several inner GMRES iterations.  
\begin{figure}[!h]
    \centering
    \includegraphics[width=1.0\textwidth, keepaspectratio]{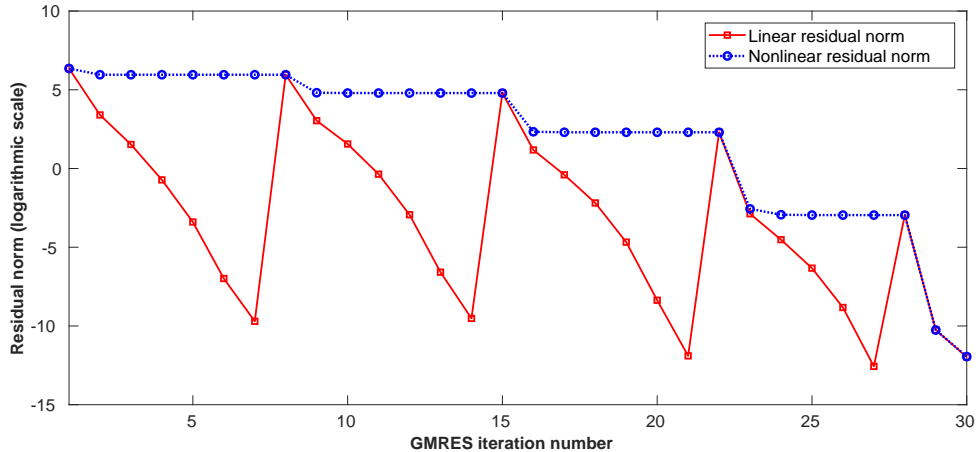}
    \caption{This figure illustrates the issue of \textit{oversolving} when solving a two-phase flow problem with 5 Newton iterations and 25 linear solver iterations.}
    \label{fig:NewtonOversolving}
\end{figure}
Figure~\ref{fig:NewtonOversolving} illustrates this phenomenon during the simulation of a coupled two-phase flow model,  wherein the red (square marker) curve shows the decrease in the linear residual over GMRES iterations for 5 Newton iterations and the blue (circle marker) curve represents the nonlinear residual norm. Away from the solution, the disagreement between the function and its local linear model is significant and this causes the stagnation of the nonlinear residual over successive linear iterations. Closer to a solution, the nonlinear residual decreases at a faster rate. Nonlinear solution algorithms can be tailored to take advantage of this phenomenon and to produce highly efficient class of solvers. 

\subsection{Inexact Newton method}
Inexact Newton method was introduced by Dembo et al. \cite{dembo1982}.
They proposed a modification of the classical Newton method by varying the degree to which the Newton equations are solved.
The inexact Newton method for discrete nonlinear equations, $\mathcal R$, is given in Algorithm~\ref{Alg:InexactNewton} with the index $\nu$ representing the Newton iteration number, $R_{tol}$ the desired residual tolerance
and $\eta_\nu$ the forcing term. 
The notation $\eta_\nu \in [0, 1[$ indicates that $\eta_\nu$ takes values between $0$ and $1$ including $0$ but excluding $1$.
\begin{algorithm}
\caption{Inexact Newton}
\label{Alg:InexactNewton}
\begin{algorithmic}
\STATE{INPUT: $u^0$}
\STATE{OUTPUT: $u^{\nu}$ such that $\norm{\mathcal R(u^{\nu})} \le R_{tol}$}
\STATE{$\nu \leftarrow 0$}
\WHILE{$\norm{\mathcal R(u^{\nu})} \le R_{tol}$}
\STATE{Compute $\eta_\nu \in [0,1[$ and $\delta^\nu$ such that}
\STATE{\begin{equation}\label{Eq:Inexact_Newton} \qquad \norm{\mathcal R \left(u^{\nu} \right) + \mathcal R'\left(u^{\nu} \right )\delta^\nu } \le \eta_\nu\norm{\mathcal R\left(u^{\nu}\right)}
\end{equation}}
\STATE{$u^{\nu+1} = u^{\nu} + \delta^\nu$}
\STATE{$\nu \leftarrow \nu + 1$}
\ENDWHILE
\RETURN $u^{\nu}$
\end{algorithmic}
\end{algorithm}
The exact \textit{Newton step}, the \textit{Newton equation} 
or the \textit{local linear model} \cite{eisenstat1996} of $R(u^{\nu+1})$
are given by the exact solution of 
 \begin{equation}\label{Eq:Newton_equation}
 \mathcal R\left(u^{\nu} \right) + \mathcal R'\left(u^{\nu} \right)\delta^\nu = 0.
 \end{equation}
  $\eta_\nu$ controls the accuracy to which the Newton equation needs to be solved. In theory, both $\eta_\nu$ and $\delta^\nu$ can be selected such that they satisfy Eq.~\eqref{Eq:Inexact_Newton}.  In most practical applications, each $\eta_\nu$ is computed first and the Newton update is computed such that Eq.~\eqref{Eq:Inexact_Newton} is satisfied. 
  In commercial reservoir simulators, $\eta_\nu$ is assumed to be a constant between $10^{-3}$ to $10^{-6}$. 
  
\subsection{Variable forcing term}  
As mentioned in Section~\ref{Sec:Intro}, solving exactly or using with very small forcing terms for the Newton equation leads to \textit{oversolving}. 
To mitigate this issue, several variable forcing terms have been proposed in the literature.
 Brown and Saad \cite{brown1990} proposed hybrid Krylov methods for the solution of nonlinear systems in which the Newton equation is solved only approximately. For the forcing term, the authors proposed an estimate $\eta_\nu = \left(\frac{1}{2}\right)^\nu$ where $\nu = 1,2,\ldots$. This choice results in a q-superlinear convergence but as it does not incorporate any information on the function $\mathcal R$, 
 we observe poor convergence rates when our test cases are highly nonlinear.
 The reader is referred to the book by Dennis and Schnabel \cite{dennis1983} for the types of convergence defined in this paper.
Eisenstat and Walker \cite{eisenstat1996}, in their seminal work, proposed two estimates of the forcing term such that the inexact Newton method resulted in little \textit{oversolving} while preserving the attractive local convergence properties of the classical Newton method. Their first choice reflects the agreement between the nonlinear function and its local linear model at the previous step and is given by
\begin{equation}\label{Eq:Choice_1_original}
    \eta_\nu = \frac{\norm{\mathcal R(u^{\nu}) - \mathcal R(u^{\nu - 1}) - \mathcal R'(u^{\nu - 1})\delta^{\nu - 1} }}{\norm{\mathcal R(u^{\nu - 1})}} \qquad \nu = 1,2,\ldots.
\end{equation}
Their second choice for the forcing term, which does not take into account the agreement between the nonlinear function and its local linear model, is based on the decay of the norm of the nonlinear residual. This is given as
\begin{equation}\label{Eq:Choice_2_original}
    \eta_\nu = \gamma \left(\frac{\norm{\mathcal{R}(u^\nu)}}{\norm{\mathcal{R}(u^{\nu - 1})}}\right)^r \qquad \nu = 1,2,\ldots,
\end{equation}
where $r \in ]1,2]$ and $\gamma \in [0,1]$. The convergence of the inexact Newton method using the second choice as the forcing term is mainly affected by the coefficient $\gamma$. A too restrictive $\gamma$ will result in \textit{oversolving} while a large value will result in an increased number of nonlinear iterations. 
Dawson et al. \cite{dawson1997} implemented these two estimates in the context of reservoir simulation for two-phase immiscible problems.
They use a sequential formulation solving separately the pressure
and the saturation equation.
They use for both systems GMRES with a two-stage Gauss-Seidel (2SGS) preconditioner.
Their results show an improvement in the number of linear iterations
but an increase in the number of nonlinear iterations.
%
Cai and Keyes \cite{cai2002} studied the convergence properties of constant forcing terms in the context of high Reynolds flows.
They proposed a preconditioned inexact Newton method based on nonlinear additive Schwarz algorithms and used GMRES combined with one-level additive Schwarz.
They used constant forcing terms between $\eta_\nu = 10^{-3}$ to $\eta_\nu = 10^{-6}$ and noticed for these specific cases what are the optimal values.
Various applications in the context of Navier-Stokes equations along with the forcing term estimates of Eqs.~\eqref{Eq:Choice_1_original} and \eqref{Eq:Choice_2_original} can be found in Hwang and Cai \cite{hwang2005}.
Their results show an improvement in the number of linear iterations
but an increase in the number of nonlinear iterations.
%
Recently, Almani \cite{almani2016} presented results 
for the inexact Newton method 
with the estimates of Eqs.~\eqref{Eq:Choice_1_original} and \eqref{Eq:Choice_2_original}
for coupled porous media flow and geomechanics problems.
The solutions were compared with a fixed forcing term of $10^{-6}$.
The linear solver used was GMRES preconditioned either by Line-SOR or AMG \cite{stuben2001a, henson2002}.
The author presented significant decrease in the accumulated linear iterations
but increase in the nonlinear iterations.
Most of the research works need to combine backtracking with inexact Newton methods in order not to deteriorate too much the global nonlinear convergence. 
An et al. \cite{an2007} proposed recently, in the context of ....,
 a strategy to predict the forcing term using the information from the actual reduction, $a_\nu$, and the predicted reduction, $p_\nu$, in the residual norm. 
$a_\nu$ and $p_\nu$ were defined in \cite{eisenstat1994} as 
\begin{eqnarray} 
    a_\nu(\delta^\nu) &=& \norm{\mathcal R(u^\nu)} - \norm{\mathcal R(u^{\nu} + \delta^\nu)} \text{ and}\\
    p_\nu(\delta^\nu) &=& \norm{\mathcal R(u^\nu)} - \norm{\mathcal R(u^{\nu}) + \mathcal R'(u^{\nu})\delta^\nu  }.
\end{eqnarray}
Defining $t_\nu = a_\nu(\delta^\nu)/p_\nu(\delta^\nu)$, the forcing term is given by 
\begin{equation}\label{Eq:Choice_AN}
    \eta_\nu = 
    \begin{cases}
    1 - 2p_1, &\text{ if } t_\nu < p_1 \\
    \eta_\nu, &\text{ if } p_1 \le t_\nu < p_2 \\
    0.8\eta_\nu, &\text{ if } p_2 \le t_\nu < p_3 \\
    0.5\eta_\nu, &\text{ if } t_\nu \ge p_3, \\
    \end{cases}
\end{equation}
where $0 < p_1 < p_2 < p_3 < 1$ are user defined parameters with $0 < p_1 < 0.5$. The authors proved q-superlinear convergence for their choice. 
Botti \cite{botti2015} proposed a new approach, for specific nonlinear functions, given by 
\begin{equation}\label{Eq:Choice_Botti}
    \eta_\nu = \frac{\norm{\mathcal R(u^{\nu}) + \mathcal R'(u^{\nu})\delta^\nu }}{\norm{\mathcal R(u^{\nu}) + \mathcal R'(u^{\nu})\delta^\nu  } + \alpha \norm{\mathcal R(u^\nu) - \mathcal R(u^{\nu} + \delta^\nu)}},
\end{equation}
where $\alpha \in ]1, 2]$, is a user defined parameter that controls the behavior of the sequence of the forcing terms. This choice is interpreted as a predictor-corrector strategy where the prediction results in satisfactory local convergence while the correction reduces the \textit{oversolving}. 
We tested these last two approaches on our problems and the results were not as interesting as the methods presented hereafter.

\subsection{Globalization}
The convergence of the Newton and the inexact Newton methods is only local. This means that for initial guesses not "sufficiently" close to the solution, the iterative process may not converge. There are several globalization techniques for the Newton method that improve the likelihood of convergence from arbitrary initial guesses. Interesting references can be found in \cite{eisenstat1994} where the authors analyze globally convergent inexact Newton algorithms such as \textit{minimum reduction method} and \textit{trust level method}.
Botti \cite{botti2015} also presents a globalization algorithm using the \textit{pseudo-transient continuation} method and successful application of the new estimate to several test cases. Globalization of the inexact Newton algorithm is not in the scope of the current work as  reservoir simulator already use various advanced theoretical and heuristic safeguarding methods such as the ones presented in \cite{ECL, jenny2009, li2014, moyner2017}. 

\subsection{Proposed work}
The aim of the present work is to develop better estimates for the forcing term, $\eta_\nu$, such that the number of iterations of the new inexact Newton method is close to the 
number of iterations of the reference inexact Newton method with constant forcing term
but with less linear iterations.
In the next section, theory related to the new forcing terms is presented,
followed by computational results from general purpose reservoir simulation. 
Test cases ranging from simultaneous water and gas injection in heterogeneous porous media to Water Alternating Gas (WAG) are presented. 
Black-oil and compositional simulation examples are presented as well as a test case of Steam-Assisted Gravity Drainage (SAGD). 
Finally, the results from a deep offshore west-african field with gas re-injection are presented.
 
\section{Method}
This section is divided into two parts. 
The proofs of the local convergence of the inexact Newton method obtained for the forcing terms proposed in this work 
and the characterization of the order of convergence for the new forcing term estimates.
These are extensions of the previous works of \cite{dembo1982, eisenstat1996}.
The proposed choices achieve fast local convergence and avoid \textit{oversolving} to an extent. All the proposed choices are scale independent and do not change if the function is multiplied by a constant.

\subsection{Preliminary development}
Certain assumptions on $\mathcal R$ are introduced along with some useful constants for the following proofs. For simplicity of notation, the iteration index, $\nu$, is dropped from the state variables and used only where necessary. 
The H\"older condition is a generalization of the Lipschitz continuity condition.
A function $\mathcal R'$ satisfies the H\"older condition at $u^*$ if there are non-negative real constants $C$ and $\alpha$ that satisfy for all $u$ in the neighborhood of $u^*$ 
\begin{equation} 
\norm{\mathcal R'(u) - \mathcal R'(u^*)} \le C \norm{u - u^*}^\alpha \label{Eq:Holder}.
\end{equation}
Usually $\alpha \in ]0, 1]$ and when $\alpha = 1$, $\mathcal R'(u)$ is Lipschitz continuous at $u^*$.
Let $u^* \in \mathbb{R}^n$ such that $\mathcal{R}(u^*) = 0$ and set 
\begin{equation}\label{Eq:Derivative_inverse}
    M \equiv \max\{\norm{\mathcal{R}'(u^*)}, \norm{\mathcal{R}'(u^*)^{-1}}\}.
\end{equation} 
We will use a region around $u^*$ that satisfies the inequality
\begin{equation}\label{Eq:condition_on_inverse}
    \norm{\mathcal R'(u)^{-1}} \le (\alpha + 1)M
\end{equation}
and we will introduce a value $\xi^*$ that satisfies
\begin{equation}\label{Eq:condition_on_radius}
    \xi^* < \frac{\alpha+1}{CM}.
\end{equation}
Let's define the ball of radius $\xi$ around $u^*$
\begin{equation}
B_\xi(u^*) = \left\{ u \in \mathbb{R}^n / \norm{u - u^*} < \xi  \right\}.
\end{equation}
Let $\xi^*$ be a sufficiently small strictly positive constant such that
\begin{itemize}
    \item $\mathcal R(u)$ is continuously differentiable and $\mathcal R'(u)$ nonsingular on $B_{\xi^*}(u^*)$,
    \item inequality \eqref{Eq:Holder} holds for $u \in B_{\xi^*}(u^*)$ with $C$ and $\alpha$ non-negative real constants,
    \item inequality \eqref{Eq:condition_on_inverse} holds for $u \in B_{\xi^*}(u^*)$,
    \item inequality \eqref{Eq:condition_on_radius} holds.
\end{itemize}
\begin{lemma}\label{Lem:Lemma1}
Given $u, v \in B_{\xi^*}(u^*)$, then 
\begin{equation}\label{Eq:Lemma_1}
    \norm{\mathcal{R}(v) - \mathcal{R}(u) - \mathcal{R}'(u)(v - u)} \le C\left ( 2\norm{u - u^*}^\alpha + \frac{\norm{v - u}^\alpha}{\alpha + 1}  \right ) \norm{v - u}.
\end{equation}
\end{lemma}
\textbf{Proof} This proof is an extension to the H\"older condition of Lemma 1.1 in \cite{eisenstat1996} valid for the Lipschitz condition ($\alpha = 1$).
A continuous equivalent of the step $u$ to $v$ is defined such that 
$u(t) \equiv u + t\delta$ for $0\le t \le 1$, where $\delta = v - u$. 
Using the fundamental theorem of calculus for line integrals
\begin{equation}
    {\mathcal{R}(v) - \mathcal{R}(u)} = {\int_0^1 \bigg [ \mathcal{R}'(u + t\delta) \bigg ] \delta dt }. 
\end{equation}
Subtracting $\mathcal{R}'(u)\delta$ from both sides in the above equation and taking the norm gives 
\begin{eqnarray}
    \norm{\mathcal{R}(v) - \mathcal{R}(u) - \mathcal{R}'(u)\delta} &=& \norm{\int_0^1 \bigg [ \mathcal{R}'(u + t\delta) - \mathcal{R}'(u) \bigg ] \delta dt } \nonumber\\
 &\le& \int_0^1 \norm{ \bigg [ \mathcal{R}'(u + t\delta) - \mathcal{R}'(u^*) + \mathcal{R}'(u^*)  - \mathcal{R}'(u) \bigg ]  dt } \norm{\delta}\nonumber\\    
 &\le& \int_0^1 \bigg [ \norm{ \mathcal{R}'(u + t\delta) - \mathcal{R}'(u^*)} + \norm{\mathcal{R}'(u^*)  - \mathcal{R}'(u)} \bigg ] dt \norm{\delta}.
\end{eqnarray}
Eq.~\eqref{Eq:Holder} substituted for the first two and the last two terms in the inequality above results in 
\begin{eqnarray}
    \norm{\mathcal{R}(v) - \mathcal{R}(u) - \mathcal{R}'(u)\delta} &\le& C \int_0^1 \bigg [ \norm{(u + t\delta) - u^*}^\alpha + \norm{u - u^*}^\alpha \bigg ] dt \norm{\delta} \nonumber\\
                                                                         &\le& C \bigg( \int_0^1  \big[ \norm{u - u^*}^\alpha + t^\alpha\norm{\delta}^\alpha \big]  + \norm{u - u^*}^\alpha \bigg) dt \norm{\delta}\nonumber\\
                                                                         &=& C \bigg( \norm{u - u^*}^\alpha + \frac{\norm{\delta}^\alpha}{\alpha + 1} + \norm{u - u^*}^\alpha \bigg) \norm{\delta}\nonumber\\
                                                                         &=&C\left( 2\norm{u - u^*}^\alpha + \frac{\norm{\delta}^\alpha}{\alpha + 1} \right) \norm{\delta}.
\end{eqnarray}    

\begin{lemma}\label{Lem:Lemma2}
    There exists a strictly positive constant $\mu > 0$, such that for any $u \in B_{\xi^*}(u^*)$, 
    \begin{equation}\label{Eq:Lemma_2}
        \frac{1}{\mu}\norm{u - u^*} \le \norm{\mathcal{R}(u)} \le \mu\norm{u - u^*}.
    \end{equation}
\end{lemma}
\textbf{Proof} 
This proof is similar to Lemma 1.2 in \cite{eisenstat1996} with the new Eq.~\eqref{Eq:Lemma_1}.
Adding and subtracting $\mathcal{R}'(u^*)(u - u^*)$ to $\mathcal{R}(u)$ we get
\begin{equation}\label{Eq:Lemma_2_p0}
\mathcal{R}(u) = \mathcal{R}(u) - \mathcal{R}'(u^*) (u - u^*) + \mathcal{R}'(u^*)(u - u^*).
\end{equation}
As $\mathcal{R}(u^*) = 0$, the term $\mathcal{R}(u^*)$ can be added to Eq.~\eqref{Eq:Lemma_2_p0} to give 
\begin{equation}\label{Eq:Lemma_2_p1}
\mathcal{R}(u) = \mathcal{R}(u) - \mathcal{R}(u^*) - \mathcal{R}'(u^*) (u - u^*) + \mathcal{R}'(u^*)(u - u^*).
\end{equation}
Taking the norm of Eq.~\eqref{Eq:Lemma_2_p1} results in
\begin{align}
\norm{\mathcal{R}(u)} &= \norm{\mathcal{R}(u) - \mathcal{R}(u^*) - \mathcal{R}'(u^*) (u - u^*) + \mathcal{R}'(u^*)(u - u^*)} \label{Eq:Lemma_2_p2} \\
                      &\le \norm{\mathcal{R}(u) - \mathcal{R}(u^*) - \mathcal{R}'(u^*) (u - u^*)} + \norm{\mathcal{R}'(u^*)(u - u^*)} \label{Eq:Lemma_2_p3}
\end{align}
Using $u=u^*$ and $v = u$ in Eq.~\eqref{Eq:Lemma_1} gives
\begin{equation}\label{Eq:Lemma_2_p4}
    \norm{\mathcal{R}(u) - \mathcal{R}(u^*) - \mathcal{R}'(u^*) (u - u^*)} \le 0.0 + C \frac{\norm{u - u^*}^{\alpha + 1}}{\alpha + 1}.
\end{equation}
Using Eqs.~\eqref{Eq:Derivative_inverse} and \eqref{Eq:Lemma_2_p4}, Eq.~\eqref{Eq:Lemma_2_p3} becomes
\begin{eqnarray}
\norm{\mathcal{R}(u)}   &\le& \norm{\mathcal{R}(u) - \mathcal{R}(u^*) - \mathcal{R}'(u^*) (u - u^*)} + \norm{\mathcal{R}'(u^*)(u - u^*)} \nonumber\\
                        &\le& C \frac{\norm{u - u^*}^{\alpha + 1}}{\alpha + 1} + M\norm{u - u^*} \nonumber\\
                        &\le&\left( C \frac{\norm{u - u^*}^{\alpha}}{\alpha + 1} + M \right) \norm{u - u^*}.
\end{eqnarray}
Now to prove that $\frac{1}{\mu}\norm{u - u^*} \le \norm{\mathcal{R}(u)}$, we write that 
\begin{equation}\norm{\mathcal{R}(u)} = \norm{- \mathcal{R}(u)},\end{equation}
we add $\mathcal{R}(u^*)=0$ and we add and subtract $\mathcal{R}'(u^*)(u - u^*)$ to get
\begin{eqnarray}
\norm{\mathcal{R}(u)}   &=& \norm{- \mathcal{R}(u) + \mathcal{R}(u^*) + \mathcal{R}'(u^*)(u - u^*)   - \mathcal{R}'(u^*) (u - u^*)} \nonumber\\
                        &\ge& - \norm{\mathcal{R}(u) - \mathcal{R}(u^*) - \mathcal{R}'(u^*) (u - u^*)} + \norm{\mathcal{R}'(u^*) (u - u^*)} \nonumber\\
                        &=& - \norm{\mathcal{R}(u) - \mathcal{R}(u^*) - \mathcal{R}'(u^*) (u - u^*)} + \norm{(\mathcal{R}'(u^*))^{-1}}^{-1}\norm{(u - u^*)}.
\end{eqnarray}
Using Eqs.~\eqref{Eq:Derivative_inverse} and \eqref{Eq:Lemma_2_p4},
\begin{eqnarray}
\norm{\mathcal{R}(u)}   &\ge& - C \frac{\norm{u - u^*}^{\alpha + 1}}{\alpha + 1} + \frac{1}{M}\norm{u - u^*}  \nonumber\\
                        &=&\left( - C \frac{\norm{u - u^*}^{\alpha}}{\alpha + 1} + \frac{1}{M}  \right) \norm{u - u^*}.
\end{eqnarray}
The lemma is then proved using 
\begin{equation}
 \mu \equiv \max \left\{ C \frac{\norm{u - u^*}^{\alpha}}{\alpha + 1} + M, \, \left(- C \frac{\norm{u - u^*}^{\alpha}}{\alpha + 1} +  \frac{1}{M} \right)^{-1} \right\}.   
\end{equation}

\begin{lemma}\label{Lem:Lemma3}
    There exists a strictly positive constant $\beta$, such that for any Newton update $\delta$ and forcing term $\eta \in [0,1[$ satisfying Eq.~\eqref{Eq:Inexact_Newton},
    we have $\forall u \in B_{\xi^*}(u^*)$, 
    \begin{equation}\label{Eq:Lemma_3}
        \norm{\mathcal{R}(u + \delta)} \le \left( \eta + \beta \norm{\mathcal{R}(u)}^\alpha \right) \norm{\mathcal{R}(u)}.
    \end{equation}
\end{lemma}

\textbf{Proof} Adding and subtracting $\mathcal{R}(u) + \mathcal{R}'(u)\delta$ to $\mathcal{R}(v)$, where $v = u + \delta$, results in
\begin{equation}\label{Eq:Lemma_3_p1}
    \mathcal{R}(v) = \mathcal{R}(v) - \mathcal{R}(u)- \mathcal{R}'(u)\delta + \mathcal{R}(u) + \mathcal{R}'(u)\delta.
\end{equation}
Taking the norms on both sides and using Lemma~\ref{Lem:Lemma1} and Eq.~\eqref{Eq:Inexact_Newton} gives
\begin{eqnarray}
    \norm{\mathcal{R}(v)} &\le& \norm{\mathcal{R}(v) - \mathcal{R}(u) - \mathcal{R}'(u)\delta} + \norm{\mathcal{R}(u) + \mathcal{R}'(u)\delta} \nonumber\\
                     &\le& C\left( 2\norm{u - u^*}^\alpha + \frac{\norm{\delta}^\alpha}{\alpha + 1} \right) \norm{\delta} + \eta \norm{\mathcal{R}(u)}  \label{Eq:Lemma_3_p2}.
\end{eqnarray}
It can be seen that 
\begin{eqnarray}
    \norm{\delta} &=& \norm{\mathcal{R}'(u)^{-1} \mathcal{R}'(u) \delta} \nonumber\\
                  &=& \norm{\mathcal{R}'(u)^{-1}}\norm{\mathcal{R}'(u) \delta} \nonumber\\
                  &\le& (\alpha + 1)M  \norm{- \mathcal{R}(u) + \mathcal{R}(u) + \mathcal{R}'(u) \delta}    \nonumber\\
                  &\le& (\alpha + 1)M\left( \norm{\mathcal{R}(u)} + \norm{\mathcal{R}(u) + \mathcal{R}'(u) \delta}\right) \nonumber\\
                  &\le& (\alpha + 1)M(1+\eta) \norm{\mathcal{R}(u)} \nonumber\\
                  &\le& 2(\alpha + 1)M \norm{\mathcal{R}(u)} \label{Eq:Lemma_3_p3},
\end{eqnarray}
which is identical to Lemma 1.3 given in \cite{eisenstat1996} if $\alpha=1$. 
Substituting Eqs.~\eqref{Eq:Lemma_3_p3} and~\eqref{Eq:Lemma_2} in Eq.~\eqref{Eq:Lemma_3_p2} results in 
\begin{eqnarray}
    \norm{\mathcal{R}(v)} &\le& C\left( 2\mu^\alpha \norm{\mathcal{R}(u)}^\alpha + \frac{(2(\alpha + 1)M)^\alpha}{\alpha + 1} \norm{\mathcal{R}(u)}^\alpha \right) 2(\alpha + 1)M\norm{\mathcal{R}(u)} + \eta \norm{\mathcal{R}(u)} \nonumber\\
                     &\le& \left[ 2(\alpha + 1)MC\left( 2\mu^\alpha + \frac{(2(\alpha + 1)M)^\alpha}{\alpha + 1} \right) \norm{\mathcal{R}(u)}^\alpha + \eta \right] \norm{\mathcal{R}(u)} \label{Eq:Lemma_3_p4},
\end{eqnarray}
which proves the lemma with 
\begin{equation}\label{Eq:Lemma_3_p5}
    \beta = 2MC\bigg( 2\mu^\alpha(\alpha + 1) + \big(2(\alpha + 1)M\big)^\alpha\bigg).
\end{equation}

Eqs.~\eqref{Eq:Lemma_3_p4} and~\eqref{Eq:Lemma_3_p5} reduce to Lemma 1.4 in \cite{eisenstat1996} for $\alpha = 1$. 
Lemmas~\ref{Lem:Lemma1}-\ref{Lem:Lemma3} will be used to deduce the rate of convergence for the particular choices of the forcing term proposed in this work.

In the following, a convergence analysis of the proposed choices is presented. For the remainder of the document we set $\xi$ such that 
\begin{equation}\label{Eq:xi_choice} 0 < \xi \le \frac{\xi^*}{1 + 2(\alpha + 1)M \mu}.
\end{equation} 
With  this inequality, Eqs.~\eqref{Eq:Lemma_2} and \eqref{Eq:Lemma_3_p3} are used to prove Lemma~\ref{Lem:Lemma_2.1} that is equivalent to Proposition 2.1 in \cite{eisenstat1996}.
\begin{lemma}\label{Lem:Lemma_2.1}
    Given $u \in B_{\xi}(u^*)$ with $\xi$ from Eq.~\eqref{Eq:xi_choice}, for all $\delta \in \mathbb{R}^n$ and $\eta \in [0,1[$ satisfying Eq.~\eqref{Eq:Inexact_Newton}, we have $u + \delta \in B_{\xi^*}(u^*)$.
\end{lemma}

\textbf{Proof} Substituting Eq.~\eqref{Eq:Lemma_2} into Eq.~\eqref{Eq:Lemma_3_p3}, 
\begin{equation}\label{Eq:Lem2.1_3}
\norm{\delta} \le 2(\alpha + 1)M \mu\norm{u - u^*}
\end{equation}
and as $u \in B_{\xi}(u^*)$, $\norm{u - u^*} \le \xi$
and as $\xi$ satisfies Eq.\eqref{Eq:xi_choice}, Eq.~\eqref{Eq:Lem2.1_3} becomes
\begin{align}
\norm{\delta} &\le 2(\alpha + 1)M \mu \, \xi \le \xi^* - \xi \label{Eq:Lem2.1_4}.
\end{align}
As $u \in B_{\xi}(u^*)$ and $\norm{\delta} \le \xi^* - \xi$, it follows that $u + \delta \in B_{\xi^*}(u^*)$. 

\subsection{Choice 1}
Given the Newton iteration count, $\nu$, and an initial forcing term, $\eta_0 \in [0, 1[$, 
our first expression for the forcing term is  
\begin{equation}\label{Eq:Choice_1}
\eta_\nu = \left( \frac{\norm{\mathcal{R}(u^\nu) - \mathcal{R}(u^{\nu-1}) - \mathcal{R}'(u^{\nu-1})\delta^{\nu-1}}}{\norm{\mathcal{R}(u^{\nu - 1})}} \right)^{p^{\nu}}, \qquad\nu = 1,2,\ldots,   
\end{equation}
where $p^0 = 1.0$ and $p^\infty = 2.0$.
Choice 1 is scale independent.
It is an extension of Eq.~\eqref{Eq:Choice_1_original} from Eisenstat and Walker \cite{eisenstat1996} where they use $p^\nu=1$.
They prove that the convergence is q-superlinear if the initial guess is sufficiently close to the solution. 
Eisenstat and Walker tried to square the forcing term ($p^\nu=2$)
but reported that it was not as successful in their experiments as the other choices.
Squaring the forcing term might generate very small values and usually results in \textit{oversolving} of the Newton equation away from the solution, i.e., when the rate of convergence is not at least r-quadratic \cite{deuflhard2004}.
On the other hand, we observe that squaring the forcing term close to the solution accelerates the local convergence, resulting in less nonlinear iterations.
To take advantage of this phenomenon, we introduce a power $p^{\nu} \in [1.0,2.0]$ that increments from $1.0$ to $2.0$ as the nonlinear iterations proceed. 
\begin{lemma}\label{Lem:Lemma4}
    Let's have the constants $C$, $\alpha$, $M$, $\xi^*$, $\mu$, $\beta$, $\xi$, ${p^{\nu}}$ retaining their definitions introduced before and set $\kappa=2(\alpha + 1)M$.
    Let $\eta_0 \in [0,1[$ and $\hat \eta$ be such that $\eta_0 < \hat \eta < 1$. Given an initial guess, $u^0 \in B_{\xi}(u^*)$, sufficiently close to the solution, and a strictly positive constant $\theta$ such that 
    \begin{itemize}
        \item $\left[ \kappa C \left( 2\mu^\alpha + \frac{\kappa^\alpha}{\alpha + 1} \right) \theta^\alpha \right]^{p^{\nu}} + \beta\theta^\alpha \le \hat\eta$, 
        \item $\eta_0 + \beta\theta^\alpha \le \hat\eta$,
        \item $\theta  < \xi/\mu$
        \item and $\norm{\mathcal{R}(u^0)} \le \theta$,
    \end{itemize}
then $u^\nu \in B_{\xi}(u^*) \subset B_{\xi^*}(u^*)$ and $\norm{\mathcal{R}(u^\nu)} \le \theta$ $\forall \nu = 0,1,2,\ldots$.
\end{lemma}

\textbf{Proof} This proof is an generalize of the one presented in \cite{eisenstat1996} with the H\"older constants and the variable exponent of the forcing term, ${p^{\nu}}$.
Lemma~\ref{Lem:Lemma3} for the first step starting from $u^0$ is written as 
\begin{align}
    \norm{\mathcal{R}(u^1)} &\le \left( \eta_0 + \beta \norm{\mathcal{R}(u^0)}^\alpha \right) \norm{\mathcal{R}(u^0)} \nonumber\\
                              &\le \left( \eta_0 + \beta \theta^\alpha \right) \norm{\mathcal{R}(u^0)} \nonumber\\
                              &\le \hat\eta \norm{\mathcal{R}(u^0)} \nonumber\\
                              &\le \theta \label{Eq:Theo_p1}
\end{align}
Since $u^0 \in B_{\xi}(u^*)$, from Lemma~\ref{Lem:Lemma_2.1}, $u^1 \in B_{\xi^*}(u^*)$.
Using the left inequality of Lemma~\ref{Lem:Lemma2} with $\norm{\mathcal{R}(u^1)} \le \theta$
results in
\begin{equation}
    \frac{1}{\mu} \norm{u^1 - u^*} \le \norm{\mathcal{R}(u^1)} \le \theta \le \frac{\xi}{\mu}
\end{equation}
or
\begin{equation}\label{Eq:Theo_p2}
    \norm{u^1 - u^*} \le \xi,
\end{equation}
which is equivalent to $u^1 \in B_{\xi}(u^*)$. 
Hence, $u^0 \in B_{\xi}(u^*)$ and $\norm{\mathcal{R}(u^0)} \le \theta$ result in
$u^1 \in B_{\xi}(u^*)$ and $\norm{\mathcal{R}(u^1)} \le \theta$.
From Lemma~\ref{Lem:Lemma_2.1}, 
$u^{\nu} \in B_{\xi}(u^*) \subset B_{\xi^*}(u^*)$
results in 
$u^{\nu+1} \in B_{\xi^*}(u^*)$.
To finalize the proof by induction of the lemma, we now need to 
prove that $\norm{\mathcal{R}(u^{\nu-1})}<\theta$ and $\norm{\mathcal{R}(u^{\nu})}<\theta$
results in $\norm{\mathcal{R}(u^{\nu+1})}<\theta$.
Using Eq.~\eqref{Eq:Choice_1} with Lemma~\ref{Lem:Lemma1} gives
\begin{eqnarray}
    \eta_\nu &=& \left( \frac{\norm{\mathcal{R}(u^\nu) - \mathcal{R}(u^{\nu-1}) - \mathcal{R}'(u^{\nu-1})\delta^{\nu-1}}}{\norm{\mathcal{R}(u^{\nu - 1})}} \right)^{p^{\nu}} \nonumber\\
             &\le& \left(\frac{C\left( 2\norm{u^{\nu-1} - u^*}^\alpha + \frac{\norm{\delta^{\nu-1}}^\alpha}{\alpha + 1} \right) \norm{\delta^{\nu-1}}}{\norm{\mathcal{R}(u^{\nu-1})}} \right)^{p^{\nu}}.
\end{eqnarray}
Using the left inequality of Lemma~\ref{Lem:Lemma2} for $\norm{u^{\nu-1} - u^*}$ and the inequality~\eqref{Eq:Lemma_3_p3} for $\norm{\delta^{\nu-1}}$ we get
\begin{eqnarray}
    \eta_\nu &\le& \left( \frac{C\left( 2\mu^\alpha \norm{\mathcal{R}(u^{\nu-1})}^\alpha + \frac{\kappa^\alpha \norm{\mathcal{R}(u^{\nu-1})}^\alpha}{\alpha + 1} \right) \kappa\norm{\mathcal{R}(u^{\nu-1})} }{\norm{\mathcal{R}(u^{\nu-1})}} \right)^{p^{\nu}} \nonumber\\
             &\le& \left[ \kappa C \left( 2\mu^\alpha + \frac{\kappa^\alpha}{\alpha + 1} \right) \norm{\mathcal{R}(u^{\nu-1})}^\alpha \right]^{p^{\nu}} \nonumber.
\end{eqnarray}
Substituting $\eta_\nu$ from the previous equation in Lemma~\ref{Lem:Lemma3} results in
\begin{eqnarray}
\norm{\mathcal{R}(u^{\nu+1})}   &\le& \left( \eta_\nu + \beta\norm{\mathcal{R}(u^\nu)}^\alpha \right) \norm{\mathcal{R}(u^\nu)} \nonumber\\
                                &\le& \left( \left[ \kappa C \left( 2\mu^\alpha + \frac{\kappa^\alpha}{\alpha + 1} \right) \norm{\mathcal{R}(u^{\nu-1})}^\alpha \right]^{p^{\nu}} + \beta\norm{\mathcal{R}(u^\nu)}^\alpha \right) \norm{\mathcal{R}(u^\nu)} \label{Eq:Theo_p4}.
\end{eqnarray}
As $\norm{\mathcal{R}(u^{\nu-1})}<\theta$ and $\norm{\mathcal{R}(u^{\nu})}<\theta$ we then have
\begin{equation}
\norm{\mathcal{R}(u^{\nu+1})}   \le \left( \left[ \kappa C \left( 2\mu^\alpha + \frac{\kappa^\alpha}{\alpha + 1} \right) \theta^\alpha \right]^{p^{\nu}} + \beta\theta^\alpha \right) \norm{\mathcal{R}(u^\nu)} \le \hat\eta \norm{\mathcal{R}(u^\nu)} \label{Eq:Theo_p4.5},
\end{equation}
which results in
\begin{equation}
\norm{\mathcal{R}(u^{\nu+1})}   \le \theta.
\end{equation}
From this result, as $\mathcal{R}(u^{0})<\theta$ and $\mathcal{R}(u^{1})<\theta$,
by direct induction $\mathcal{R}(u^{\nu})<\theta$ $\forall \nu \ge 0$.
Now, we have 
$u^{\nu} \in B_{\xi}(u^*) \subset B_{\xi^*}(u^*)$,
$u^{\nu+1} \in B_{\xi^*}(u^*)$
and 
$\mathcal{R}(u^{\nu})<\theta$ $\forall \nu \ge 0$.
Following the same steps used to derive Eq.~\ref{Eq:Theo_p2}, for any general $\nu$, we have
\begin{equation}
    \frac{1}{\mu} \norm{u^{\nu+1} - u^*} \le \norm{\mathcal{R}(u^{\nu+1})} \le \theta \le \frac{\xi}{\mu}
\end{equation}
or
\begin{equation}
    \norm{u^{\nu+1} - u^*} \le \xi,
\end{equation}
which is equivalent to $u^{\nu+1} \in B_{\xi}(u^*)$. Hence, $u^\nu \in B_{\xi}(u^*) \subset B_{\xi^*}(u^*)$ and $\norm{\mathcal{R}(u^\nu)} \le \theta$ for all $\nu = 1,2,\ldots$. 

\begin{theorem}\label{Theo:1}
Given the conditions in Lemma~\ref{Lem:Lemma4}, the iterates, $u^\nu, \nu = 1,2,\ldots$, produced by Algorithm~\ref{Alg:InexactNewton} with the forcing term given by Eq.~\ref{Eq:Choice_1} remain within $B_{\xi^*}(u^*)$ and converges to $u^*$ such that 
\begin{equation}\label{Eq:Theo_1}
    \norm{u^{\nu+1} - u^*} \le \Phi \max\left\{ \norm{u^{\nu-1} - u^*}^{\alpha . p^{\nu}}, \norm{u^\nu - u^*}^\alpha  \right\} \norm{u^\nu - u^*},
\end{equation}
where 
\begin{equation}\label{Eq:Theo_2}
     \Phi = \left( \left[ \kappa C \left( 2\mu^\alpha + \frac{\kappa^\alpha}{\alpha + 1} \right)  \mu^{\alpha} \right]^{p^{\nu}} + \beta \mu^{\alpha} \right)  \mu.
\end{equation}
\end{theorem}

\textbf{Proof} Lemma~\ref{Lem:Lemma4} proves that $u^\nu \in B_{\xi}(u^*) \subset B_{\xi^*}(u^*)$ for all $\nu = 0,1,2,\ldots$. From Eq.~\ref{Eq:Theo_p4.5}, for all $\nu = 0,1,\ldots$, $$\norm{\mathcal{R}(u^{\nu+1})} \le \hat \eta \norm{\mathcal{R}(u^{\nu})},$$ and as $\hat \eta < 1.0$, $\norm{\mathcal{R}(u^{\nu})}$ converges to zero, which in turn from the left inequality of Lemma~\ref{Lem:Lemma2}, means that $u^\nu$ converges to $u^*$. 
Using the two inequalities of Lemma~\ref{Lem:Lemma2} in Eq.~\eqref{Eq:Theo_p4} and the notation $\epsilon_{\nu}=\norm{u^{\nu} - u^*}$, results in  
\begin{equation}\label{Eq:Theo_p5}
\frac{1}{\mu}\epsilon_{\nu+1}
\le 
\left( \left[ \kappa C \left( 2\mu^\alpha + \frac{\kappa^\alpha}{\alpha + 1} \right)  \mu^{\alpha} \epsilon_{\nu-1}^\alpha \right]^{p^{\nu}} + \beta \mu^{\alpha} \epsilon_{\nu}^\alpha \right)  \mu \epsilon_{\nu}.
\end{equation}
Rearranging Eq.~\eqref{Eq:Theo_p5} gives
\begin{equation}\label{Eq:Theo_p6}
\epsilon_{\nu+1} \le \Phi \max \left( [\epsilon_{\nu-1}^{\alpha}]^{p^{\nu}}, \epsilon_{\nu}^\alpha \right) \epsilon_{\nu},
\end{equation}
where 
\begin{equation}\label{Eq:Theo_p7}
   \Phi = \left( \left[ \kappa C \left( 2\mu^\alpha + \frac{\kappa^\alpha}{\alpha + 1} \right)  \mu^{\alpha} \right]^{p^{\nu}} + \beta \mu^{\alpha} \right)  \mu. 
\end{equation}
Eqs.~\eqref{Eq:Theo_p6} and~\eqref{Eq:Theo_p7} show r-quadratic convergence for a Lipschitz continuous function if $\alpha = 1$ and ${p^{\nu}} = 2$. Eqs.~\eqref{Eq:Theo_p6} reduces to Eq. 2.3 in \cite{eisenstat1996} for $\alpha = 1$ and ${p^{\nu}} = 1$.  

Choise 1 given by Eq.~\ref{Eq:Choice_1} considers the agreement between the nonlinear function and its local linear model. Some cases might result in a really good agreement far away from the solution which makes the forcing term significantly small. This situation will result in significant \textit{oversolving} of the Newton equations. To counter this situation, a second choice is presented that does not take the agreement between the function and its local linear model into account. It takes into consideration the ratio of the norm of the previous residual to the current one. This choice is described in detail in the next subsection.  

\subsection{Choice 2}
Choice 2 that is proposed in this work is inspired from practical experience on different problem settings. 
Application of Eq.~\eqref{Eq:Choice_2_original} to reservoir simulation problems resulted in more linear iterations than Eq.~\eqref{Eq:Choice_1_original} but resulted generally in less nonlinear iterations. 
Our proposed modified Choice 2 is given by 
\begin{equation}\label{Eq:Choice_2}
    \eta_\nu = \phi(\nu) \left(\frac{\norm{\mathcal{R}(u^\nu)}}{\norm{\mathcal{R}(u^{\nu - 1})}}\right)^r, \qquad \nu = 1,2,\ldots,
\end{equation}
where $r \in ]1,2]$ and $\phi$ is a function that is monotonically decaying with the Newton iteration count $\nu$. 
Several options for this function, $\phi(\nu)$, are presented in Section~\ref{Sec:Param}, along with its effects on the performance of Algorithm~\ref{Alg:InexactNewton} in Section~\ref{Sec:Results}.
Using the forcing term given by Eq.~\eqref{Eq:Choice_2} in Algorithm~\ref{Alg:InexactNewton} results in a highly desirable rate of convergence with a reduction of the cumulative number of Newton iterations. 
The proof of the effective rate of convergence is identical to the one presented in Eisenstat and Walker \cite{eisenstat1996} other than the fact that instead of a constant coefficient, $\phi_0$, an exponentially decaying function is considered. Only the case of $\phi_0 < 1$ is considered in this work keeping in view the practical applications. The theorem and its proof are presented in Appendix~\ref{App:Choice2}. 

\section{Computational results}
Performance of Algorithm~\ref{Alg:InexactNewton} applied to reservoir simulation is explored in this section. An example of \textit{oversolving} of the Newton equations in two-phase immiscible flow is shown in Subsection~\ref{Sec:Oversolving}. Various options for the variable coefficients in Eqs.~\eqref{Eq:Choice_1} and~\eqref{Eq:Choice_2} are presented in Subsection~\ref{Sec:Param}. The test problems range from heterogeneous two-phase flow models to multiphase multicomponent models containing steam as well as a real example from a west-African field. 

\subsection{Oversolving example}
\label{Sec:Oversolving}
\begin{figure}[!ht]
     \begin{center}
        \subfigure[Inexact Newton method with a fixed relative tolerance of 1e-4]{
            \includegraphics[width=0.97\textwidth]{Figures/newton_oversolving.eps}
            \label{fig:Oversolving_newt} 
        }
        \subfigure[Inexact Newton method with variable forcing term]{
            \includegraphics[width=0.97\textwidth]{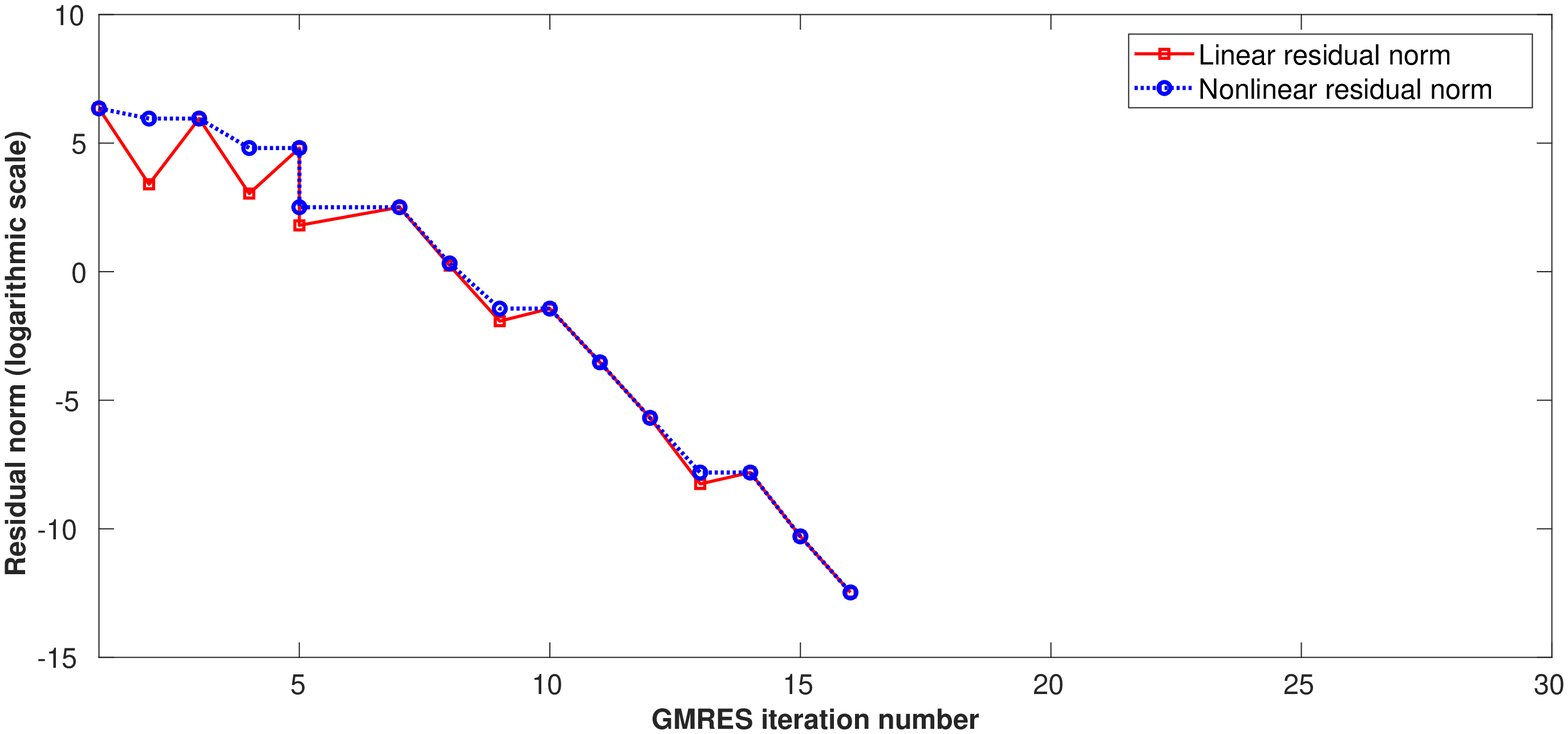}
            \label{fig:Oversolving_inex} 
        }
  	\end{center}
    \caption{ Illustration of effect of the inexact Newton method on oversolving. \ref{fig:Oversolving_newt} shows the oversolving resulting from an inexact Newton method with a fixed relative tolerance when applied to a two-phase flow problem. The blue line (circle markers) is the nonlinear residual over each GMRES iteration. The red line (square markers) is the linear residual obtained from the GMRES routine. \ref{fig:Oversolving_inex} is obtained by the application of the inexact Newton method with the forcing term Choice 1 of Eq.~\eqref{Eq:Choice_1_original} onto the same problem. For this specific time-step this choice results in little to no oversolving.}
    \label{fig:Oversolving}
\end{figure}

An example of \textit{oversolving} is presented in Figure~\ref{fig:Oversolving}. The results are obtained by solving the saturation part of a two-phase sequential implicit simulator using ILU(0)-GMRES as the linear solver. In this figure, the blue (circle marker) curve is obtained by evaluating the nonlinear residual at each linear step while the red (square marker) curve is obtained by evaluating the linear residual at each GMRES iteration. Figure~\ref{fig:Oversolving_newt} shows the behavior of inexact Newton method with a fixed relative tolerance of $10^{-4}$ where the decrease in the linear residual doesn't necessarily result in the reduction of the nonlinear residual. On the other hand, Figure~\ref{fig:Oversolving_inex} shows the effect of variable forcing term using Eq.~\eqref{Eq:Choice_1_original} applied to adapt the convergence tolerance of the GMRES solver. This method results in little to no oversolving which implies a reduced number of linear iterations. 

\subsection{Variable parameters, $p^\nu$ and $\phi$}
\label{Sec:Param}
For the original expressions of the forcing terms of Eqs.~\eqref{Eq:Choice_1_original} and~\eqref{Eq:Choice_2_original}, we use the values of $p^\nu = 1.0$, $\gamma = 0.5$ and $r = 1.618$.
We explain here the various ways of varying the parameter $p^{\nu}$ from Choice~1 given in Eq.~\eqref{Eq:Choice_1}. 
Initially, away from the solution, this parameter is taken to be equal to unity
and as the iterations proceed, $p^{\nu}$ is monotonically increased to the value 
$2.0$ to ensure quadratic convergence. 
Three different rates of changes are tested to decrease $p^\nu$, 
these choices are given by
\begin{align}
    &\text{InEx 1Steep: } p^\nu = \min \left(2.0, 2.0-\frac{2.5}{\nu}\exp(-\nu)\right) \label{Eq:C1_op1}\\
    &\text{InEx 1exp: } p^\nu = \min \left(2.0, 2.0-\exp(1 - \nu^{0.7})\right) \label{Eq:C1_op2}\\
    &\text{InEx 1cub: } p^\nu =  \min \left(2.0, \frac{\nu^3}{250.0} + \frac{\nu^2}{250.0} + \frac{\nu}{250.0} + 1\right) \label{Eq:C1_op3},
\end{align}
where $\nu$ is the Newton iteration number. 
The labels used in front of the formulas will be used in the remainder of the text to specify the choice being used. 

\begin{figure}[!ht]
     \begin{center}
        \subfigure[Convergence rate, $p$]{
            \includegraphics[width=0.47\textwidth]{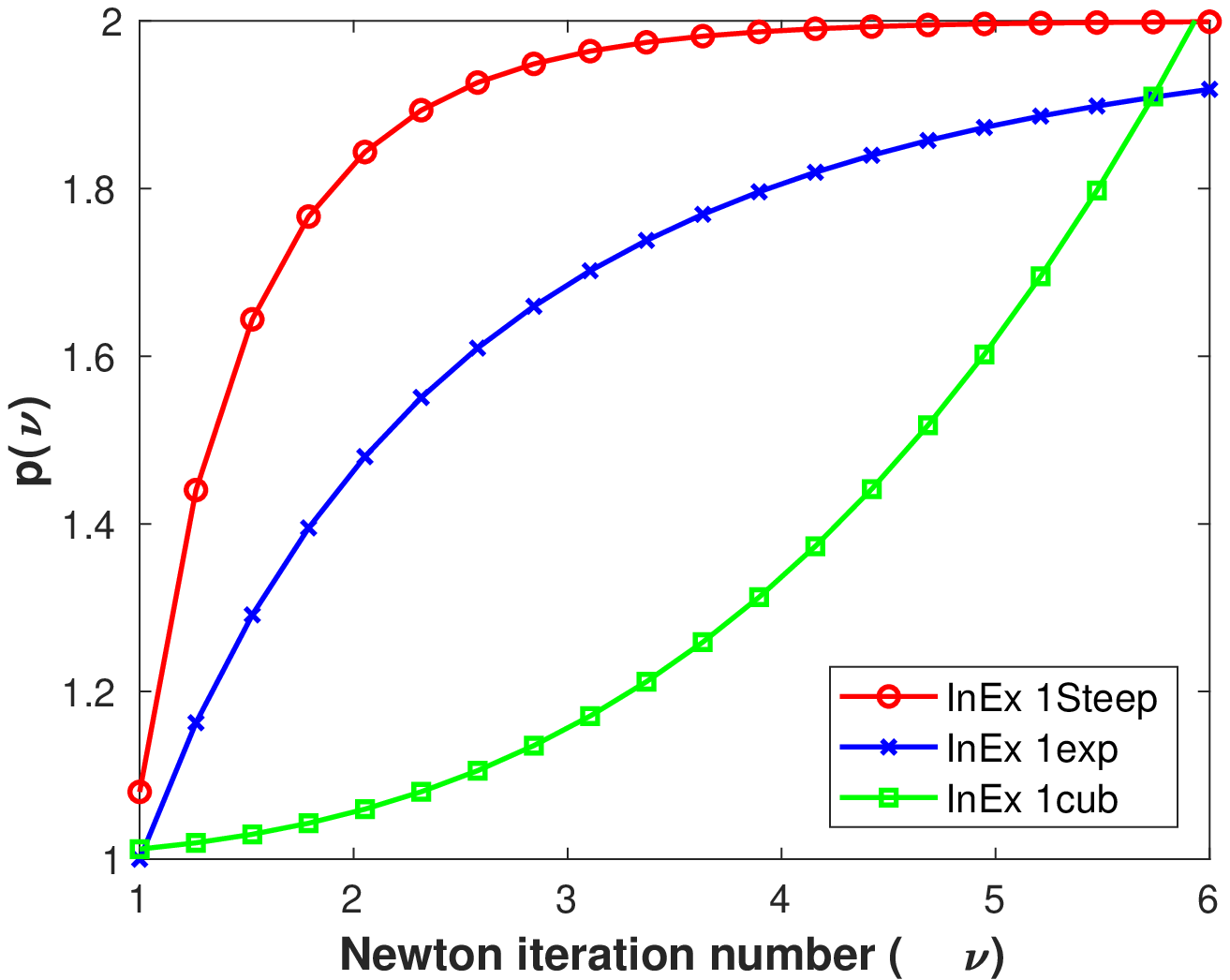}
            \label{fig:Choice_1_p} 
        }
        \subfigure[Variable coefficient, $\phi$]{
            \includegraphics[width=0.47\textwidth]{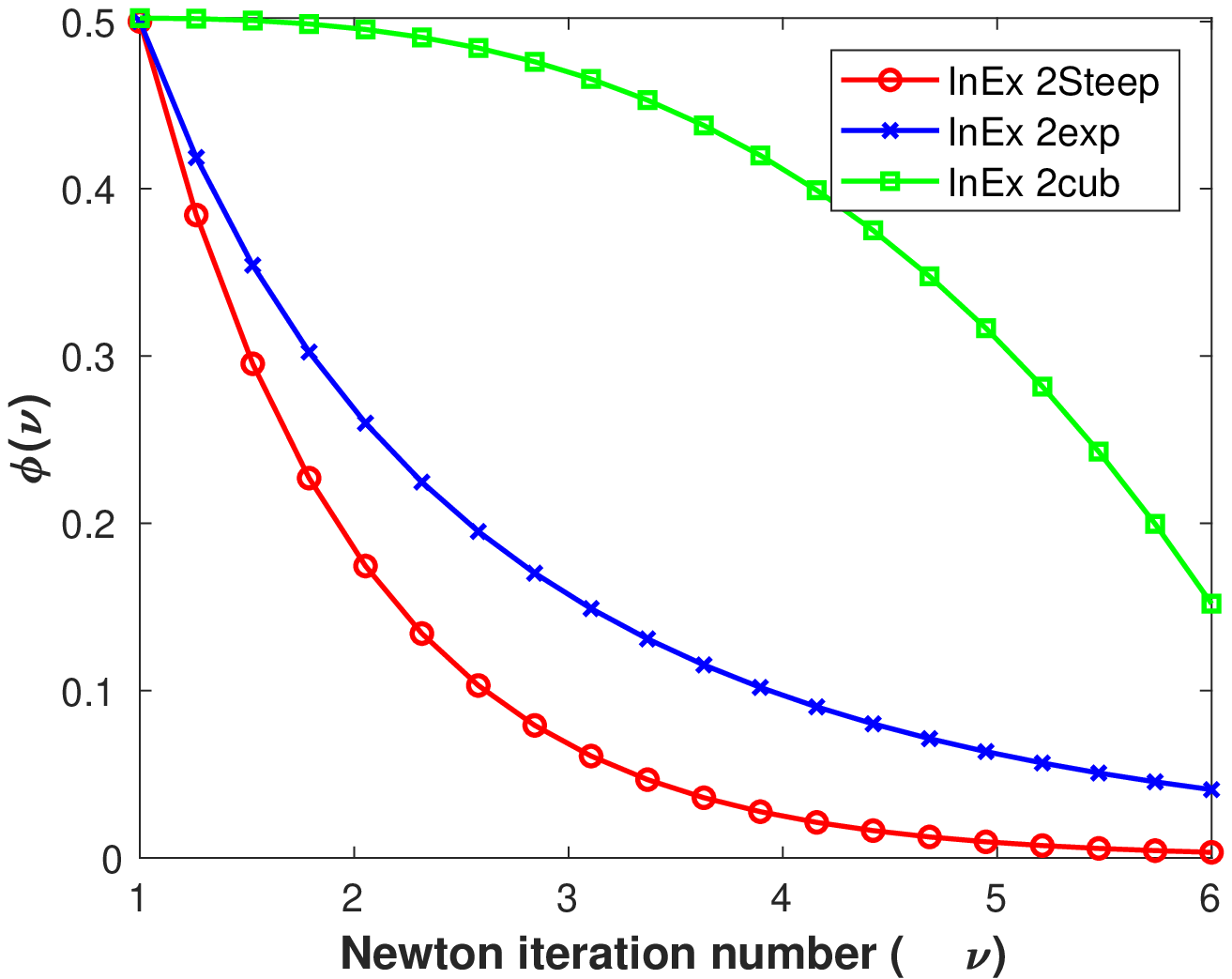}
            \label{fig:Choice_2_phi} 
        }
  	\end{center}
    \caption{Figure shows the behaviour of the variable coefficients, $p$ and $\phi$, in Eqs.~\ref{Eq:Choice_1} and~\ref{Eq:Choice_2}}
    \label{fig:VariableCoefficients}
\end{figure}

Choice~2 in Eq.~\eqref{Eq:Choice_2} has one parameter, 
$r \in ]1,2]$ usually set to $r = 1.618$ and $\phi$ a monotonically decaying function 
with the Newton iteration count $\nu$. 
The following formula list the choices explored for $\phi$: 
\begin{align}
    &\text{InEx 2Steep: } \phi(\nu) = \max \left(\epsilon_{0}, \phi_0 \exp(1-\nu) \right) \label{Eq:C2_op1}\\
    &\text{InEx 2exp: } \phi(\nu) = \max \left(\epsilon_{0}, \phi_0\exp(1-\nu^{0.7}) \right) \label{Eq:C2_op2}\\
    &\text{InEx 2cub: } \phi(\nu) =  \max \left(\epsilon_{0}, \phi_0 \left(- \frac{\nu^3}{250.0} + \frac{\nu^2}{250.0} + \frac{\nu}{250.0} + 1\right) \right) \label{Eq:C2_op3}, 
\end{align}
where $\epsilon_{0}$ is the minimum tolerance for the linear solver. 
The behaviour of the above choices is shown in Figure~\ref{fig:VariableCoefficients}.   
The red (circle marker), blue (cross) and green (square) curves
show the results for both $p^\nu$ and $\phi(\nu)$
with the the steep exponential choice, the gradual exponential and the cubic strategies, respectively.

\subsection{Application to reservoir simulation}
\label{Sec:Results}
The modified choices for the forcing term presented in the previous sections along with the variable coefficient choices are implemented 
in the company in-house research reservoir simulator \cite{Moncorge:2012,Jaure:2014}.
These test cases include black-oil, compositional and thermal models in various geological settings. 

Tables~\ref{Tab:Newt_compare} to~\ref{Tab:InEx_2_compare} summarize the results for the different forcing term choices. 
Table~\ref{Tab:Newt_compare} shows the cumulative linear and nonlinear iterations for fixed forcing term between $\eta = 10^{-6}$ to $\eta = 10^{-1}$.
Table~\ref{Tab:InEx_1_compare} shows the cumulative linear and nonlinear iterations for
the original and new Choice 1 of Eqs.~\eqref{Eq:Choice_1_original} and Eqs.~\eqref{Eq:Choice_1}.
Table~\ref{Tab:InEx_2_compare} shows the cumulative linear and nonlinear iterations for
the original and new Choice 2 of Eqs.~\eqref{Eq:Choice_2_original} and Eqs.~\eqref{Eq:Choice_2}.
We compare the results in Tables~\ref{Tab:InEx_1_compare} and \ref{Tab:InEx_2_compare}
with the fixed forcing term Newton method with $\eta = 10^{-4}$.
InEx 1 and InEx 2 are the original methods 
of Eqs.~\eqref{Eq:Choice_1_original} and \eqref{Eq:Choice_2_original}
with $\gamma = 0.5$ and $r = 1.618$
while 
InEx 1Steep, InEx 2Steep,
InEx 1Exp, InEx 2Exp,
InEx 1Cub and InEx 2Cub
are the new methods
of Eqs.~\eqref{Eq:Choice_1} and \eqref{Eq:Choice_2}
with $p^\nu$ given by Eqs.~\eqref{Eq:C1_op1}, \eqref{Eq:C1_op2}, \eqref{Eq:C1_op3},
$\phi(\nu)$ given by Eqs.~\eqref{Eq:C2_op1}, \eqref{Eq:C2_op2}, \eqref{Eq:C2_op3},
$r = 1.618$ and $\epsilon_0=10^{-6}$.
The first number in each column is the cumulative number of GMRES iterations using the CPR-AMG preconditioner and the second number inside the brackets gives the total number of Newton iterations including time-step cuts. 

Comparing the results for the seven test cases, we observe that inexact Newton method with variable forcing term results in little oversolving and thus result in less number of linear iterations. Even though the methods take considerably less linear iterations, the cumulative number of nonlinear iterations tend to increase. In all the cases tested in this work, except the thermal steam problem, Choice 1 of Eq.\eqref{Eq:Choice_1_original} resulted in the least number of linear iterations. On the other hand, Choice 2 of Eq.\eqref{Eq:Choice_2_original} show better nonlinear convergence with less cumulative Newton iterations. 
The modified choices, given by Eqs.~\eqref{Eq:Choice_1} and~\eqref{Eq:Choice_2}, invariably resulted in better nonlinear convergence to their original counterparts. 
The InEx 2Steep method has the best convergence properties along with resulting between 30\% to 60\% savings in the linear iterations. 
For the thermal steam problem, only InEx 2Steep produces acceptable results as the other methods take considerably higher number of nonlinear iterations.

\begin{table}[!h]
\centering
\begin{tabular}{|c|l|l|l|l|l|}
\hline
Cases   & 1e-6                & 1e-4                 & 1e-3                & 1e-2        & 1e-1\\ 
\hline
IPAM1  &13365(\textbf{1612}) &9040(\textbf{1734})   &6180(\textbf{1618})  &    4276 (\textbf{1679})   & 2981 (\textbf{1943})          \\ 
SIAM5  &48077(\textbf{6643}) &34835(\textbf{7455})  &23532(\textbf{6614}) &    16878 (\textbf{7105})   &  21593 (\textbf{8039})        \\
IPAM2  &2117(\textbf{157})   &1404(\textbf{157})    &1044(\textbf{157})   &    642  (\textbf{143})     &  377 (\textbf{140})            \\
ANTBO  &13201(\textbf{1664}) &8172(\textbf{1664})   &5899(\textbf{1663})  & 4113 (\textbf{1674})     & 2844 (\textbf{1856})     \\ 
ANTC7  &27689(\textbf{3378}) &17409(\textbf{3378})  &12729(\textbf{3400}) & 8928 (\textbf{3398})    & 5805 (\textbf{3623})            \\
SAGD   &12121(\textbf{1010}) &8822(\textbf{976})    &7899(\textbf{1099})  &    6852(\textbf{1287})   &  5733(\textbf{1736})        \\
BIGC    &322791(\textbf{7641}) &283848(\textbf{6670})    &150717(\textbf{6149}) &  91527 (\textbf{6048})   & 111868 (\textbf{7216})            \\
\hline
\end{tabular}
\caption{Cumulative linear and nonlinear iterations for fixed forcing term Newton methods. The numbers in normal font show the cumulative linear iterations while the numbers in bracket and bold font give the total nonlinear iterations required. }
\label{Tab:Newt_compare}
\end{table}
\begin{table}[!h]
\centering
\begin{tabular}{|c|l|l|l|l|l|}
\hline
Cases   & 1e-4 (ref.)            & InEx 1                   & InEx 1Steep                  & InEx 1Exp                 & InEx 1Cub \\
\hline
IPAM1  & 9040(\textbf {1734})   & 2760(\textbf {1780})  & 3548(\textbf{1727})    &3446(\textbf{1834})   &2927(\textbf{1789})  \\
SIAM5  & 34835(\textbf {7455})  & 17606(\textbf {8920}) & 20568(\textbf{8198})   &14955(\textbf{7601})  &15455(\textbf{8609}) \\
IPAM2  & 1404(\textbf {157})    & 346(\textbf {158})    & 396(\textbf{142})      &399(\textbf{148})     &384(\textbf{148})    \\
ANTBO  & 8172(\textbf{1664})    & 2997(\textbf{1853})   & 3829(\textbf{1808})    &3621(\textbf{1826})   &3148(\textbf{1847})  \\
ANTC7  & 17409(\textbf {3378})  & 6025(\textbf {3719})  & 8025(\textbf {3585})   &7510(\textbf{3602})   &6443(\textbf{3618})  \\
SAGD   & 8822(\textbf {976})    & 7696(\textbf {3442})  & 6610(\textbf {2449})   &6631(\textbf{2496})   &6850(\textbf{2665})  \\   
BIGC    & 283848(\textbf {6670}) & 90148(\textbf {7763}) & 53092(\textbf {6869})  &-                     &54356(\textbf{7158})  \\   
\hline
\end{tabular}
\caption{Cumulative linear and nonlinear iterations for 
fixed forcing term Newton method $10^{-4}$,
forcing term computed with Eq.~\eqref{Eq:Choice_1_original}
and 
forcing term computed with our new estimate Eq.~\eqref{Eq:Choice_1}
with $p^\nu$ given from Eqs.~\eqref{Eq:C1_op1},~\eqref{Eq:C1_op2} and~\eqref{Eq:C1_op3}.}
\label{Tab:InEx_1_compare}
\end{table}
\begin{table}[!h]
\centering
\begin{tabular}{|c|l|l|l|l|l|}
\hline
Cases   & 1e-4 (ref.)            & InEx 2                   & InEx 2Steep                  & InEx 2Exp                 & InEx 2Cub \\
\hline
IPAM1  & 9040(\textbf {1734})   &3794(\textbf {1739}) &4068(\textbf {1580})   &4613(\textbf {1870})   &4020(\textbf {1733})  \\
SIAM5  & 34835(\textbf {7455})  &21035(\textbf {8286})&21926(\textbf {7399})  &18508(\textbf {7764})  &18559(\textbf {7984}) \\
IPAM2  & 1404(\textbf {157})    &421(\textbf {160})   &564(\textbf {138})     &492(\textbf {139})     &449(\textbf {142})    \\
ANTBO  & 8172(\textbf{1664})    &3889(\textbf {1788}) &4277(\textbf {1778})   &4117(\textbf {1777})   &3905(\textbf {1788})  \\
ANTC7  & 17409(\textbf {3378})  &8363(\textbf {3552}) &9183(\textbf {3503})   &8888(\textbf {3523})   &8495(\textbf {3552})  \\
SAGD   & 8822(\textbf {976})    &6439(\textbf {2321}) &5286(\textbf {1213})   &5871(\textbf {1459})   &5821(\textbf {1502})  \\   
BIGC    & 283848(\textbf {6670}) &48537(\textbf {7298}) &129632(\textbf {6248})   &134479(\textbf {6668})   &123617(\textbf {6464})  \\   
\hline
\end{tabular}
\caption{Cumulative linear and nonlinear iterations for 
fixed forcing term Newton method $10^{-4}$,
forcing term computed with Eq.~\eqref{Eq:Choice_2_original} with $\gamma=0.5$ and $r=1.618$
and 
forcing term computed with our new estimate Eq.~\eqref{Eq:Choice_2}
with $\phi(\nu)$ given from Eqs.~\eqref{Eq:C2_op1},~\eqref{Eq:C2_op2} and~\eqref{Eq:C2_op3}
and $r=1.618$.}
\label{Tab:InEx_2_compare}
\end{table}

\subsubsection{Test case 1: WAG (IPAM1)}
This case shows the simulation of Water Alternating Gas (WAG) injection 
2D Case-1 from \cite{MoncorgeSFI:2018}.
Table~\ref{Tab:Newt_compare} compares the performance of the various inexact Newton method with fixed forcing terms. 
The stricter the convergence criteria of the linear solver, the larger the number of cumulative linear iterations. Table~\ref{Tab:Newt_compare} shows that a fixed tolerance of $10^{-6}$ takes a much larger number of GMRES iterations. This is quite a trivial observation but the difference between these curves indicate the degree of \textit{oversolving} in the solution of the Newton equations. The larger the discrepancy between these curves, the larger the \textit{oversolving}.

 As the convergence tolerance gets loser, the number of cumulative GMRES iterations decrease. An inverse trend is expected for the case of Newton iterations but ad-hoc time-integration strategies based on the number of linear iterations could result in minor discrepancies. Looking at the number of nonlinear iterations for this case, $10^{-2}$ takes less iterations than $10^{-4}$, which is due to smaller number of time-steps for the former. Comparison of the variable $\eta_\nu$ with the fixed values is done only for $10^{-3}$ and $10^{-4}$ as they are the most practically used criteria. Table~\ref{Tab:InEx_1_compare} shows that InEx 1 results in 55\% and 60\% reduction in the GMRES iterations compared to $10^{-3}$ and $10^{-4}$, respectively. On the other hand, it takes 10\% and 2.6\% more nonlinear iterations. Using InEx 1Steep of Eq.~\eqref{Eq:C1_op1} results in a slight increase in the linear iterations than the InEx 1 choice but shows better convergence rate. InEx 1Cub, which uses the cubic increase for $p$ in Eq.~\eqref{Eq:Choice_2}, shows very similar results to InEx 1. Table~\ref{Tab:InEx_2_compare} shows the best results. 
 InEx 2Steep of Eq.~\eqref{Eq:C2_op1} reduces the number of linear iterations by 35\% compared to $10^{-3}$ but at the same time also reduces the number of nonlinear iterations by 2.4\%. 
 It is also noted that the slow exponential change (option 2) for both Choice 1 and Choice 2 was not successful and resulted in a lot of nonlinear iterations.   

\subsubsection{Test case 2: WAG (SIAM5)}
This case has the same scenario than IPAM1 but is using a more complex compositional fluid.
Table~\ref{Tab:Newt_compare} shows the numbers obtained when a fixed relative tolerance is given to the GMRES solver. These results show expected trends except for the case with $\eta_\nu = 0.1$, where at one step the number of linear iterations increases steeply. This is because the tolerance is too big for the Newton method so converge efficiently. 
Out of the methods in the Choice 1 class (Table~\ref{Tab:InEx_1_compare}), 
InEx 1exp performed the best reducing the linear iterations by 57\% and increasing cumulative inexact Newton iterations by only 15\%, compared to the case with fixed $\eta_\nu = 10^{-3}$. 
Similar results were obtained for the class of estimates given by Choice 2 (Table~\ref{Tab:InEx_2_compare}) where the slow exponential decay, Eq.~\eqref{Eq:C2_op2}, of $\phi(\nu)$ resulted in a decrease in the number of GMRES iterations while keeping the nonlinear iterations in check. In both these cases, Choice 1 with steep exponential change given by Eq.~\eqref{Eq:C1_op1} performed better than the case with $\eta_\nu = 10^{-4}$, while producing comparable cumulative GMRES iterations with respect to the less strict tolerance of $10^{-3}$. It is interesting to note that the case with fixed $\eta_\nu = 10^{-3}$ produced the most efficient results for this particular problem but that might not be the true for other settings as it will be seen in the following results. Hence, its important to develop an efficient and robust criteria that works without the need for much user expertise.

\subsubsection{Test case 3: Water/gas injection (IPAM2)}
This case shows the water gas injection of 2D Case-2 from \cite{MoncorgeSFI:2018}.
In this case, water and gas were injected simultaneously at different locations in the reservoir model. Table~\ref{Tab:Newt_compare} shows the cumulative number of GMRES and inexact Newton iterations taken for the simulation of this problem. The tighter the linear convergence tolerance, the larger the number of linear iterations. Choosing $\eta_\nu = 0.1$ results in 5.6 times less cumulative linear iterations than $\eta_\nu = 10^{-6}$. This shows a great degree of \textit{oversolving}. For such a case, a variable forcing term should theoretically produce very efficient results. As can be seen from Table~\ref{Tab:InEx_1_compare}, InEx 1Steep takes two-fold less linear iterations than a fixed choice of $\eta_\nu = 10^{-3}$, which is almost four-fold less than the tightest convergence tolerance. Along with reducing the linear iteration count, it results in 13\% less nonlinear iterations. Similar results are seen for Choice 2 and its variants given in Table~\ref{Tab:InEx_2_compare}. The nonlinear convergence of all the variants of Choice 2 is comparable while the linear iteration count is reduced by at least two-fold compared to the fixed choice of $10^{-3}$.

\subsubsection{Test case 4: Anticline BO}
This case shows the 3D case from \cite{MoncorgeSFI:2018}.
A synthetic anticlinal model is setup which is tested against a black-oil fluid type. The degree of \textit{oversolving} can be seen in Table~\ref{Tab:Newt_compare}. The choice $\eta_\nu = 0.1$ obviously results in the least number of linear iterations but in this case not an ideal nonlinear convergence rate. Every other fixed method results in a slightly improved number of nonlinear iterations. 

Variants of Choice 1 (Table~\ref{Tab:InEx_1_compare}) 54\% or more reduction in the number of GMRES iterations compared to $\eta_\nu = 10^{-3}$, while resulting in a convergence rate quite similar to the option with fixed $\eta_\nu = 0.1$. A slightly better nonlinear convergence rate is obtained by using Choice 2 and its variants (Table~\ref{Tab:InEx_2_compare}), while still resulting in around 40\% reduction in the the number of linear iterations. All the variable forcing term choices presented here produce acceptable results.

\subsubsection{Test case 5: Anticline C7}
This case shows the 3D case from \cite{MoncorgeECMOR:2018,Moncorge:2019}.
The tests are done on the anticlinal model used in the previous case with a compositional fluid with 7 thermodynamic components. 
Table~\ref{Tab:Newt_compare} is obtained by running the inexact Newton algorithm with several fixed forcing terms. $\eta_\nu = 0.1$ results in 4.5 times less linear iterations and comparable Newton convergence rate. Table~\ref{Tab:InEx_1_compare} and~\ref{Tab:InEx_2_compare} show the results for the application of Choice 1 and Choice 2 along with their variants to the inexact Newton algorithm, respectively. The original Choice 1 (InEx 1), produced the least number of linear iterations but a slightly larger number of nonlinear iterations. On the other hand, for the case shown in Table~\ref{Tab:InEx_2_compare}, InEx 2 and its modifications result in similar number of iterations for both linear and the nonlinear loops. The gain in the number of linear iterations, compared to the case with $\eta_\nu = 10^{-3}$, is around two-fold.

\subsubsection{Test case 6: Thermal steam simulation: SAGD}
The final test case describes a thermal problem containing steam. 
This is a Steam Assisted Gravity Drainage (SAGD) 2D test case with one point injector injecting steam and one point producer.
This is a very complicated problem as far as the nonlinear convergence is considered. 
A very hot steam chamber is developing during the simulation.
Compared to other cases, the degree of \textit{oversolving} of the Newton equations is not too pronounced with a fixed forcing term of $10^{-6}$ takes only twice the number of linear iterations than $0.1$ (see Table~\ref{Tab:Newt_compare}). 
For comparison, in the previous cases, this ratio was somewhere close to 5.0. 
On the other hand, the lose tolerance choice almost doubles the number of inexact Newton steps over the entire simulation. Therefore, in this case, rather than tackling oversolving very aggressively, a balanced strategy is required to keep a check on the nonlinear convergence rate as well. 
The original Choice 1 of Eq.\eqref{Eq:Choice_1_original} does not reduce the number of linear iterations compared to the choice with $\eta_\nu = 10^{-3}$. In fact, this algorithm compares very poorly to the fixed cases (Table~\ref{Tab:InEx_1_compare}) and results in about 3.5 times the number of nonlinear iterations. The variants given by Eqs.~\eqref{Eq:C1_op1} to~\ref{Eq:C1_op3}, result in similar performance degradation. Forcing term choice that takes the local linear model into consideration is unsuccessful for this thermal model. 

Choice 2 of Eq.~\eqref{Eq:Choice_2} along with its variants given by Eqs.~\eqref{Eq:C2_op1} to~\ref{Eq:C2_op3}, resulted in a reduced number of linear iterations along with a fairly acceptable nonlinear convergence rate. The original Choice 2 of Eq.~\eqref{Eq:Choice_2_original} reduces the number of GMRES iterations but a fixed coefficient in Eq.~\eqref{Eq:Choice_2} worsens the nonlinear convergence rate. The choice represented by Eq.~\eqref{Eq:C2_op1} resulted in the best gains, reducing the number of linear iterations by 1.5 folds and resulting in only 10\% extra inexact Newton iterations.

\section{Real world application: Deep offshore field in west-Africa}
This test case is a real deep offshore west African field with gas re-injection.
It contains 104,000 cells, 7 hydrocarbon components, 
Peng-Robinson corrected equation of state,
7 producers controlled by group liquid target rate and
4 injectors injecting water and produced gas.
Tables~\ref{Tab:InEx_1_compare} and~\ref{Tab:InEx_2_compare} show the cumulative linear and nonlinear iterations taken to simulate this field for $20$ years. Table~\ref{Tab:Newt_compare} compares the results for different fixed forcing term values varying from $10^{-6}$ to $10^{-3}$. 
As expected, the stricter the linear tolerance, the more the number of linear iterations. But unlike the other cases presented in this work, the number of nonlinear iterations does not follow similar trend. In fact, the choice with $\eta_\nu = 10^{-6}$ takes the maximum number of Newton iterations and this is due to more time-step cuts. $\eta_\nu = 10^{-3}$ produce the best results for the fixed strategy. 

Tables~\ref{Tab:InEx_1_compare} and~\ref{Tab:InEx_2_compare} compare the results obtained for the original forcing term estimates given by Eqs.~\eqref{Eq:Choice_1_original} and~\eqref{Eq:Choice_2_original} with the modifications presented in this work. The original estimates are able to reduce the number of linear iterations, i.e. reduce the \textit{oversolving} of the Newton equations, but result in the highest number of nonlinear iterations. It will be seen further in the section that this results in slower simulation. All the modifications presented in this work, given by Eqs.~\eqref{Eq:Choice_1} and~\eqref{Eq:Choice_2}, result in less nonlinear iterations than the original estimates. The best results are obtained by using Eq.~\eqref{Eq:Choice_1} with the exponent given by Eq.~\eqref{Eq:C1_op1}. This combination resulted in the least number of linear iterations along with keeping a check on the increase of the number of Newton iterations. 

\begin{figure}[!ht]
     \begin{center}
            \includegraphics[width=0.97\textwidth]{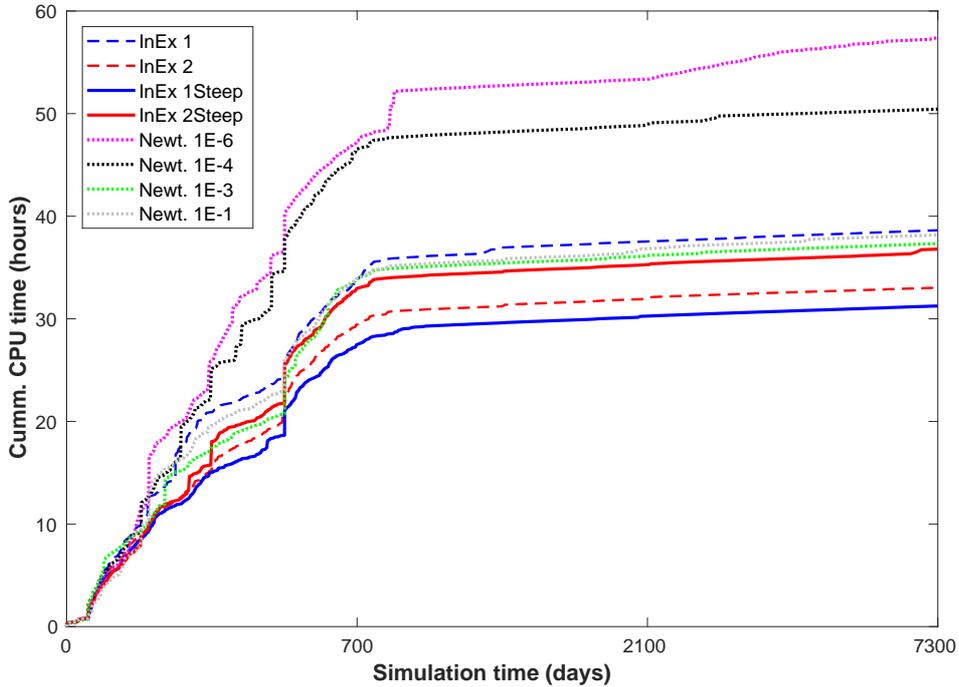}
  	\end{center}
    \caption{CPU time for inexact Newton method with different forcing term choices for a west-African deep offshore field.}
    \label{fig:time_bigC}
\end{figure}

Figure~\ref{fig:time_bigC} shows the cumulative CPU time for the simulation of the deep offshore field. Inexact Newton with fixed forcing term values of $10^{-6}$ and $10^{-4}$ takes approximately $57$ and $50$ hours to complete the simulation of $20$ years. This is a direct consequence of Table~\ref{Tab:Newt_compare}, where the choice with $\eta_\nu = 10^{-4}$ takes less linear as well as nonlinear iterations. As described before, the stricter strategy takes more nonlinear iterations due to more time-step cuts (failures). Because it is a stricter strategy, it is expected to take more linear iterations as the GMRES solver needs to converge to a smaller tolerance. Although this trend continues for $\eta_\nu = 10^{-3}$, it is not the case for a lose tolerance of $\eta_\nu = 0.1$. $\eta_\nu = 0.1$ takes less linear iterations but almost $1000$ more Newton iterations than $\eta_\nu = 10^{-3}$. This results in an increase in the computational time. As illustrated from this case, relaxing the linear tolerance does not necessarily result in faster simulation. An adaptive strategy is required. These observations stress the need for a choice that strikes a balance between the cumulative linear and nonlinear iterations. Inexact Choice 1 of Eq.~\eqref{Eq:Choice_1_original} proves to be slower in terms of computational time than the choice with fixed $\eta_\nu = 10^{-3}$. 
Even though Eq.~\eqref{Eq:Choice_1_original} takes less linear iterations, it results in a considerable increase in the nonlinear iterations. The modification presented in Eq.~\eqref{Eq:C1_op1} for Choice 1 Eq.~\eqref{Eq:Choice_1} results in the fastest simulation taking around $30$ hours. This is almost twice as fast as a fixed $\eta_\nu = 10^{-6}$ and $20\%$ faster than the best fixed choice of $\eta_\nu = 10^{-3}$. 
On the other hand, the original forcing term described by Eq.~\eqref{Eq:Choice_2_original} resulted in a faster simulation compared to its modification, given by Eq.~\eqref{Eq:Choice_2} with $\phi(\nu)$ given by Eq.~\eqref{Eq:C2_op1}, because the latter resulted in a significant increase in the cumulative linear iterations even though it required the lesser nonlinear iterations to converge. 
To conclude, forcing term choice given by Eq.~\eqref{Eq:Choice_1} resulted in an optimum balance between the linear and nonlinear iterations.  

\section{Conclusions}
In these work, we propose new forcing term choices compared to the ones presented by Eisenstat and Walker \cite{eisenstat1996}.
Proofs of convergence of the inexact Newton method for H\"older continuous functions as well as local convergence rates have been derived.
Theoretically, the proposed choices should be adapted based on the distance from the solution, but, for the application to this work, the choices are adapted based on the number of nonlinear iterations.
Our new estimates have been tested on several cases with black-oil and compositional fluids as well as a thermal model with steam.
They have also been tested on a west-African field model.
Our new estimates combined with the existing globalization techniques used in commercial reservoir simulators improve the overall convergence rate of the inexact Newton method.
We recommend using the choice given by Eq.~\eqref{Eq:Choice_2} with $\phi$ given by Eq.~\eqref{Eq:C2_op1} that provided consistent results in all the test cases.
The future work to improve the performance of the method will be to adapt the forcing term with the information from the physics of the underlying solution. 

\section*{Acknowledgements}
The authors would like to thank TOTAL management for permission to publish this work.

\appendix
\section{Convergence proof for $\eta_\nu$ given by Eq.~\eqref{Eq:Choice_2}} 
\label{App:Choice2}

The following theorem is similar to the one presented in Theorem 2.3 in Eisenstat and Walker \cite{eisenstat1996} except the variable nature of $\phi(\nu)$ given by Eq.~\eqref{Eq:Choice_2}. 
We present this theorem and proof for the purpose of completeness. 
\begin{theorem}\label{Theo:2}
     Let's have the constants $\alpha$, $\xi^*$, $\mu$, $\beta$ and $\xi$ retaining
     their definitions introduced before,
     $r$ the exponent described in Eq.~\ref{Eq:Choice_2}, 
     $\phi(\nu)$ described in Eq.~\ref{Eq:Choice_2} with $\phi(0)=\phi_0<1$ decreasing,
     $\eta_0 \in [0,1[$
     and 
     an initial guess $u^0 \in B_{\xi}(u^*)$.
     Let's define a sufficiently small strictly positive constant $\theta$ such that 
    \begin{itemize}
        \item $\eta_0 + \beta \theta^\alpha \le \eta_0^{1/r}$,
        \item $\theta < \xi/\mu$
        \item and $\norm{\mathcal{R}(u^0)} \le \theta$,
\end{itemize}
    then the following iterates, $u^\nu, \nu = 1,2,\ldots$, produced by Algorithm~\ref{Alg:InexactNewton} with the forcing term given by Eq.~\eqref{Eq:Choice_2} remains within $B_{\xi^*}(u^*)$ and converges to $u^*$ with q-order r. 
\end{theorem}

\textbf{Proof} Following the statement of Theorem~\ref{Theo:2}, Lemma~\ref{Lem:Lemma_2.1} gives $u^1 \in B_{\xi^*}(u^*)$. Similar to the first part of the proof for Lemma~\ref{Lem:Lemma4} (Eqs.~\ref{Eq:Theo_p1} to~\ref{Eq:Theo_p2}), we prove that $u^1 \in B_{\xi}(u^*)$. 

As an inductive hypothesis suppose that, for some $\nu \ge 0$, we have $u^{\nu-1} \in B_{\xi}(u^*)$, $u^\nu \in B_{\xi}(u^*)$, $\norm{\mathcal{R}(u^{\nu-1})} \le \theta$, $\norm{\mathcal{R}(u^{\nu})} \le \theta$ and $\eta_\nu < \eta_0$.
Let's prove that we have $u^{\nu+1} \in B_{\xi}(u^*)$
and $\norm{\mathcal{R}(u^{\nu+1})} \le \theta$.
Lemma~\ref{Lem:Lemma_2.1} gives $u^{\nu+1} \in B_{\xi^*}(u^*)$. Lemma~\ref{Lem:Lemma3} applied at iterate $\nu$ gives
\begin{eqnarray}
    \norm{\mathcal{R}(u^{\nu+1})} &\le& \left( \eta_\nu + \beta\norm{\mathcal{R}(u^\nu)}^\alpha \right) \norm{\mathcal{R}(u^\nu)} \nonumber\\
                                    &\le& \left( \eta_0 + \beta \theta^\alpha \right) \norm{\mathcal{R}(u^\nu)} \nonumber \\
                                    &\le& \eta_0^{1/r} \norm{\mathcal{R}(u^\nu)} \label{Eq:Theo_2_p0}
\end{eqnarray}
and as $\eta_0\le1.0$ and $\norm{\mathcal{R}(u^\nu)}\le\theta$ we have
\begin{equation}
    \norm{\mathcal{R}(u^{\nu+1})} \le \theta \label{Eq:Theo_2_p1}.
\end{equation}
Using Lemma \ref{Lem:Lemma2}, we get
\begin{equation}
    \frac{1}{\mu}\norm{u^{\nu+1}-u^*}\le \norm{\mathcal{R}(u^{\nu+1})} \le \theta \le \frac{\xi}{\mu} \label{Eq:Theo_2_p1_1},
\end{equation}
then giving $u^{\nu+1} \in B_{\xi}(u^*)$. 
Therefore by direct induction, we have $\forall \nu$, $u^{\nu} \in B_{\xi}(u^*) \subset B_{\xi^*}(u^*)$ and $\norm{\mathcal{R}(u^\nu)}\le\theta$.
Furthermore, using the inequality~\eqref{Eq:Theo_2_p0} in Eq.~\eqref{Eq:Choice_2}, we get 
\begin{eqnarray}
    \eta_{\nu + 1} &=& \phi(\nu + 1) \left(\frac{\norm{\mathcal{R}(u^{\nu + 1})}}{\norm{\mathcal{R}(u^{\nu})}}\right)^r \nonumber \\
                   &\le& \phi(\nu + 1) \eta_0 \le \phi_0 \eta_0 \le \eta_0 \label{Eq:ForcingTermDecay}. 
\end{eqnarray}
As $\eta_0\le1.0$, Eq.~\eqref{Eq:Theo_2_p0} proves that $\norm{\mathcal{R}(u^{\nu})}$ converges to zero which in turn,
by using Lemma~\ref{Lem:Lemma2}, proves that $u^\nu$ converges to $u^*$.
Substituting Eq.~\eqref{Eq:Choice_2} into Eq.~\eqref{Eq:Lemma_3} gives 
\begin{equation}\label{Eq:Theo_2_p2}
        \norm{\mathcal{R}(u^{\nu+1})} \le \left( \phi(\nu) \left(\frac{\norm{\mathcal{R}(u^\nu)}}{\norm{\mathcal{R}(u^{\nu - 1})}}\right)^r + \beta\norm{\mathcal{R}(u^\nu)}^\alpha \right) \norm{\mathcal{R}(u^\nu)}.
\end{equation}
Since for all $\nu$, $\phi(\nu) \le \phi_0$, 
\begin{equation}\label{Eq:Theo_2_p3}
        \norm{\mathcal{R}(u^{\nu+1})} \le \left( \phi_0 \left(\frac{\norm{\mathcal{R}(u^\nu)}}{\norm{\mathcal{R}(u^{\nu - 1})}}\right)^r + \beta\norm{\mathcal{R}(u^\nu)}^\alpha \right) \norm{\mathcal{R}(u^\nu)}.
\end{equation}
Dividing Eq.~\eqref{Eq:Theo_2_p3} by $\norm{\mathcal{R}(u^\nu)}^r$ and setting $\rho_\nu = \norm{\mathcal{R}(u^\nu)}/\norm{\mathcal{R}(u^{\nu-1})}^r$ gives
\begin{equation}\label{Eq:Theo_2_p4}
    \rho_{\nu+1} \le \phi_0 \rho_\nu + \beta \norm{\mathcal{R}(u^\nu)}^{\alpha + 1 - r}.    
\end{equation}
Eq.~\eqref{Eq:Theo_2_p4} is also giving 
\begin{equation}\label{Eq:Theo_2_p5}
    \rho_{\nu+1} \le (\phi_0)^{\nu} \rho_1 +  \beta \sum_{k=0}^{\nu-1} (\phi_0)^k \norm{\mathcal{R}(u^{\nu-k})}^{\alpha + 1 - r}.
\end{equation}
As $\forall \nu$, $\norm{\mathcal{R}(u^{\nu})} \le \norm{\mathcal{R}(u^{0})}$ and $\phi_0 < 1$, Eq.~\eqref{Eq:Theo_2_p5} becomes
\begin{eqnarray}
    \rho_{\nu+1} &\le& (\phi_0)^{\nu} \rho_1 +  \beta \norm{\mathcal{R}(u^{0})}^{\alpha + 1 - r} \sum_{k=0}^{\nu-1} (\phi_0)^k \nonumber\\
                 &\le& (\phi_0)^{\nu} \rho_1 +  \beta \norm{\mathcal{R}(u^{0})}^{\alpha + 1 - r} \frac{1-(\phi_0)^{\nu}}{1 - \phi_0} \nonumber\\
                 &\le& \rho_1 +  \beta \norm{\mathcal{R}(u^{0})}^{\alpha + 1 - r} \frac{1}{1 - \phi_0}. \label{Eq:Theo_2_p6}
\end{eqnarray}
It means that $\rho_{\nu}$ is uniformly bounded, then $\norm{\mathcal{R}(u^{\nu})}$ converges to $0$ with q-order r and it follows from Lemma~\ref{Lem:Lemma2} that $u^\nu$ converges to $u^*$ with q-order r as well.

\bibliographystyle{plain}
\bibliography{References}

\end{document}